\documentclass[12pt]{article}
\usepackage[dvips]{graphicx}
\usepackage{amsmath,amsfonts,amssymb,amsthm,color}

\usepackage[numbers]{natbib}

\newtheorem{theorem}{Theorem}
\newtheorem{lemma}{Lemma}
\newtheorem{corollary}{Corollary}

\newtheorem{proposition}{Proposition}

\newtheorem{definition}{Definition}
\newtheorem{example}{Example}

\usepackage{bm}
\bmdefine{\Bt}{t}
\bmdefine{\BX}{X}
\bmdefine{\BY}{Y}
\bmdefine{\BZ}{Z}
\bmdefine{\BB}{B}
\bmdefine{\BM}{M}
\bmdefine{\BD}{D}
\bmdefine{\Bi}{i}
\bmdefine{\Bj}{j}
\bmdefine{\Bx}{x}
\bmdefine{\By}{y}
\bmdefine{\Bz}{z}
\bmdefine{\Bv}{v}
\bmdefine{\Bw}{w}
\bmdefine{\Bn}{n}
\bmdefine{\Ba}{a}
\bmdefine{\Bb}{b}
\bmdefine{\Bc}{c}
\bmdefine{\Be}{e}
\bmdefine{\Bu}{u}
\bmdefine{\Bp}{p}
\bmdefine{\Bzero}{0}
\bmdefine{\Bone}{1}

\newcommand{\Q}{{\mathbb Q}}

\newcommand{\Z}{{\mathbb Z}}

\newcommand{\cI}{{\cal I}}

\newcommand{\cP}{{\cal P}}
\newcommand{\cQ}{{\cal Q}}

\newcommand{\facet}{\mathop{{\rm facet}}}
\newcommand{\fst}{{\rm fst}}

\newcommand{\rowdim}{\nu}
\newcommand{\coldim}{p}

\newcommand{\wreath}{\mathop{\rm wr}}

\title{Perturbation method for determining the group of invariance of hierarchical models}
\author{
Tomonari Sei\\
Graduate School of Information Science and Technology\\
University of Tokyo, Japan \\ 
Satoshi Aoki\\
Department of Mathematics and Computer Science\\
Kagoshima University, Japan\\
and\\
Akimichi Takemura\\
Graduate School of Information Science and Technology\\
University of Tokyo, Japan
}

\date{January 2009}

\begin{document}
\maketitle

\vspace*{10mm}\noindent
\underline{Corresponding author}\\
Tomonari Sei

\vspace*{2mm}\noindent
\underline{Address}\\
Department of Mathematical Informatics,\\
Graduate School of Information Science and Technology,\\
The University of Tokyo,\\
7-3-1, Hongo, Bunkyo-ku, Tokyo, 113-8656, Japan.

\vspace*{2mm}\noindent
\underline{Phone \& Fax}\\
+81-3-5841-6942

\vspace*{2mm}\noindent
\underline{E-mail}\\
sei@stat.t.u-tokyo.ac.jp

\newpage

\begin{abstract}
  We propose a perturbation method for determining the (largest) group
  of invariance of a toric ideal defined in \cite{largest}.  In the
  perturbation method, we investigate how a generic element in the
  row space of the configuration defining a toric ideal is mapped by a
  permutation of the indeterminates.  Compared to the proof in
  \cite{largest} which was based on stabilizers of a subset of
  indeterminates, the perturbation method gives a much simpler proof
  of the group of invariance.  In particular, we determine the group
  of invariance for a general hierarchical model of contingency tables
  in statistics, under the assumption that the numbers of the levels
  of the factors are generic.  We prove that it is a wreath product
  indexed by a poset related to the 
  intersection poset of the maximal interaction effects of the model.
%  We also prove various basic facts on hierarchical models.
\end{abstract}

\vspace*{4mm}

\vspace*{4mm}
\noindent {\it Math Subject Classification (2000):} 62H17, 05E20.

\newpage

\setlength{\baselineskip}{24pt}

\section{Introduction}
Since the introduction of the notion of Markov basis by
\cite{diaconis-sturmfels}, toric ideals associated with various
statistical models have been intensively investigated by both
statisticians and algebraists.  In particular, statistical models for
contingency tables have been rich sources for new developments (e.g.\
\cite{aoki-takemura-2003anz}, \cite{dobra-sullivant},
\cite{ohsugi-hibi-ramanujan}).  The most important statistical model
for contingency tables is the hierarchical model
(e.g.\ \cite{lauritzen1996}), which describes interactions of factors in
terms an abstract simplicial complex.  The configuration and the toric
ideal associated with a hierarchical model is highly symmetric.
Therefore it is of considerable interest to determine the (largest)
group of invariance of a general hierarchical model.  The group of
invariance is the set of permutations of the cells of contingency
tables (or the indeterminates of a polynomial ring) which leaves the
model (or, equivalently, the kernel of the configuration, or
the row space of the configuration) invariant.
Once the group of invariance is determined, a Markov
basis (or equivalently a system of binomial generators of the toric
ideal) can be very concisely described
(\cite{aoki-takemura-2008aism,largest}, \cite{sizetwo}) by a list of
representative elements from the orbits of the group.  Without the
consideration of symmetry, Markov bases for statistical problems tend
to be very large (e.g.\ \cite{hemmecke-malkin}).

Given a particular statistical model it is often easy to guess 
a candidate group, under which the model is clearly invariant.
However as shown in \cite{largest} it is often difficult to prove that
it is the largest group of invariance, i.e., every permutation outside
the group does not leave the model invariant.  In this paper we propose a
perturbation method to determine the group of invariance.  In this
method, we look at a generic element of the model and check if a 
permutation maps the element to another element in the model.  The candidate
group is shown to be the largest group of invariance, if every
permutation which maps a sufficiently generic element of the model
into the model is necessarily an element of the candidate group.
In order to show the effectiveness of this approach,
we determine the group of invariance for 
a general hierarchical model of
contingency tables, under the assumption that the
numbers of the levels of the factors are generic.
We prove that the group of invariance is a wreath product
indexed by a poset related to the 
intersection poset of the maximal interaction effects 
of the hierarchical model.

Here we give a simple illustrative example.
Consider a hierarchical model for
four-factor contingency tables with numbers of levels $I_1$, $I_2$, $I_3$ and $I_4$
and the set of facets $\{\{1,2\},\{2,3\},\{3,4\}\}$
(see Section~\ref{sec:preliminaries} for details of notation and terminology).
% The row space of the configuration is defined by
% $L=\{x\in\mathbb{Q}^{I_1\times I_2\times I_3\times I_4}\mid
%   x(i_1,i_2,i_3,i_4)=a(i_1,i_2)+b(i_2,i_3)+c(i_3,i_4)\}$,
% where $a$, $b$ and $c$ run over all two-way tables.
Our main theorems (Theorem~\ref{thm:main} and Theorem~\ref{thm:wreath}) state that under a weak regularity condition
the group of invariance of this model is
generated by the permutations of $i_2$, the permutations of $i_3$,
the permutations of $i_1$ depending on $i_2$
and the permutations of $i_4$ depending on $i_3$.
This group is strictly larger than the direct-product group 
$S_{I_1}\times S_{I_2} \times S_{I_3} \times S_{I_4}$ of 
permutations of levels for each factor.
Other examples are given in Section~\ref{sec:wreath}.
In particular, we present an example
such that  the number of orbits 
in the minimal Markov basis is smaller 
under the action of the group of invariance 
than under the action of the direct-product group (see Example~\ref{example:Markov}).

In our proof we need to establish some basic facts on hierarchical
models, which are not found in the existing statistical literature.
These facts are of independent interest and we present them in Section
\ref{sec:hierarchical-facts}.

The organization of the paper is as follows.  In Section
\ref{sec:preliminaries} we give preliminaries and present a
perturbation lemma.  In Section \ref{sec:group-of-invariance-hierarchical} 
we state our main theorem,
which expresses the group of invariance of a hierarchical model
as an intersection of wreath products of symmetric groups. 
In Section \ref{sec:hierarchical-facts}
we establish some basic facts on hierarchical models and 
in Section \ref{sec:proof}
we give a proof of the main theorem.
In Section 
\ref{sec:wreath} we rewrite the group of invariance  as a wreath product
indexed by a poset related to the 
intersection poset of the maximal interaction effects
of the hierarchical model.
We conclude the paper by some discussions in Section \ref{sec:discussions}.

\section{Preliminaries and a perturbation lemma}
\label{sec:preliminaries}

In this section we summarize preliminary facts on hierarchical models
for contingency tables, define the group of invariance  and present a
perturbation lemma, which is essential for our proofs.  We mainly
follow the notation and terminology of \cite{lauritzen1996}.

\subsection{Preliminaries on hierarchical models for contingency tables}

A hierarchical model for $m$-factor contingency tables with numbers
of levels $I_1, \dots, I_m$ 
is specified by
an abstract simplicial complex.  Let $\Delta$ be an abstract
simplicial complex (\cite[Section 2.1]{kozlov}) of subsets of a finite
set $\{1,\dots,m\}=[m]$ of ``factors''.
% Any simplex $a\in \cA$ 
% corresponds to interaction effects among the factors $\delta\in a$.
We denote the set of maximal simplices of $\Delta$ by $\facet(\Delta)=
\{ D_1, \dots, D_K\}$. Maximal simplices are called maximal
interaction effects of the model.
For each factor $j\in [m]$, the set of ``levels''of $j$ is
denoted by $\cI_j=\{1,\dots, I_j\}=[I_j]$, where
$I_j \ge 2$. The direct product of the set of levels
$\cI=\cI_1 \times \dots\times \cI_m$ is the set of ``cells'' and its element
$\Bi=(i_1, \dots, i_m)$ is a cell.  A contingency table 
$x=(x(\Bi))_{\Bi\in \cI}$ is a vector of nonnegative integers 
indexed by the cells. The number $x(\Bi)$ is the frequency of the cell $\Bi$.
In this paper, the symbol $A\subset B$ means that $A$ is a subset of $B$.
If $A$ is a proper subset of $B$, then we write $A\subsetneq B$.

For a subset $D\subset [m]$ of factors, let $\cI_D=\prod_{j\in D}
\cI_j$.
A subvector of indices $\Bi_D=(i_j)_{j\in D} \in \cI_D$ is
called ``a marginal cell''.  When a particular cell $\Bi=(i_1, \dots,i_m)$ is given, 
$\Bi_D$ is regarded as a subvector of $\Bi$, i.e., the projection of $\Bi$ onto the
coordinates in $D$.  For a contingency table $x$, its $D$-marginal table 
$x^+_D=(x^+(\Bi_D))_{\Bi_D \in \cI_D}$ is defined by 
\[
x^+(\Bi_D)=\sum_{\Bj\in \cI, \ \Bj_D=\Bi_D} x(\Bj).
\]
Similar notation is used even when $x(\Bi)$ is not necessarily a nonnegative integer.

Fix $\cI$ and a hierarchical model $\Delta$ with 
$\facet(\Delta)=\{D_1, \dots, D_K\}$.  
Write $\rowdim=\sum_{k=1}^K |\cI_{D_k}|$ and $\coldim=|\cI|$.
For each $\Bi=(i_1,\dots,i_m)\in \cI$ consider 
the following vector (cf.\ \cite{ohsugi-hibi-ramanujan}) 
\[
{\bf e}^{(1)}(\Bi_{D_1}) \oplus 
{\bf e}^{(2)}(\Bi_{D_2}) \oplus 
\dots \oplus 
{\bf e}^{(K)}(\Bi_{D_K}) \ 
\in \; \Z^\rowdim
\]
where ${\bf e}^{(k)}(\Bi_{D_k})$ is a unit coordinate
vector of dimension $|\cI_{D_k}|$ with 1 at the position
$\Bi_{D_k}$ and 0 everywhere else. The configuration $A_\Delta$
for $\Delta$ is the set of $\coldim$ vectors
\[
A_\Delta= \big\{ 
{\bf e}^{(1)}(\Bi_{D_1}) \oplus 
\dots \oplus 
{\bf e}^{(K)}(\Bi_{D_K})
\big\}_{\Bi\in \cI} .
\]
In this paper we regard $A_\Delta$ as a 
$\rowdim \times \coldim$ integral matrix
representing a linear map
from $\Q^\coldim$ to $\Q^\rowdim$. 
The matrix $A_\Delta$ can also be expressed by Kronecker products
of identity matrices and vectors  consisting of 1's
(\cite[Section 2.1]{takemura-aoki-2004aism}).
We also assume that 
the domain $\Q^\coldim$ is equipped with the standard inner product
and we identify $\Q^\coldim$ with its dual space by the standard inner product.

Let $\{ u_\Bi\}_{\Bi\in \cI}$ be the set of indeterminates indexed by
the cells and let $\{t^{(1)}_{\Bi_{D_1}} \}_{\Bi_{D_1}\in \cI_{D_1}}
\cup \dots \cup \{t^{(K)}_{\Bi_{D_K}} \}_{\Bi_{D_K}\in \cI_{D_K}}$
denote the set of indeterminates indexed by the rows of
$A_\Delta$. The toric ideal $I_{A_\Delta}$ is the kernel of the
polynomial homomorphism $\pi_\Delta$ defined by $\pi_\Delta(u_\Bi)=
t^{(1)}_{\Bi_{D_1}}\times \dots \times t^{(K)}_{\Bi_{D_K}} $.
The structure of the toric ideal is much more difficult than the
kernel of matrix $A_\Delta$.  However we will define the invariance
property  of $I_{A_\Delta}$ in terms of the invariance property of 
the kernel of $A_\Delta$.

As we discuss in Section  \ref{subsec:group-of-invariance} we are interested in
the kernel % $\ker A_\Delta$ 
of $A_\Delta$ and the linear space 
spanned by the rows of $A_\Delta$. 
In the following we denote the kernel of $A_\Delta$ and
the linear space spanned by the rows of $A_\Delta$ by 
$\ker A_\Delta$ and $r(A_\Delta)$, respectively.  Note that
$\ker A_\Delta$ and $r(A_\Delta)$ are orthogonal complements to each
other: $r(A_\Delta)=(\ker A_\Delta)^\perp$.

In statistical theory, $r(A_{\Delta})$
corresponds to a log-linear model of cell probabilities, where the
canonical parameter vector of the exponential family is specified to
lie in the linear space $r(A_\Delta)$.
We use the single term ``model'' for $\Delta$, $r(A_{\Delta})$ and
$\ker A_{\Delta}$
because they correspond to each other.

% For the configuration $A_\Delta$ associated with the hierarchical model
% $\Delta$, the 
The explicit form of $\ker A_\Delta$ and $r(A_\Delta)$ are well known 
in the literature on contingency tables
(e.g.\ \cite{lauritzen1996}). The set $\ker A_\Delta$  is written as
\begin{equation}
\label{eq:kern1}
\ker A_\Delta = \{ y \mid y^+(\Bi_D)=0, \forall \Bi_D \in {\cal I}_D,
\forall D\in \facet(\Delta) \}.
\end{equation}
For $D\subset [m]$, let $\theta_D: \cI_D \rightarrow \Q$
denote a function defined on the set of marginal cells $\cI_D$.
Then extend the domain of $\theta_D$ to $\cI$ by 
$\theta_D(\Bi)=\theta_D(\Bi_D)$. We call $\theta_D$ a function (or a table)
depending only on the marginal cell $\Bi_D$.
Let
$L_D=\{\theta_D\}\subset \Q^\coldim$ denote the linear space of
these tables. Then
\begin{equation}
\label{eq:row-space}
r(A_\Delta)= \sum_{D\in\facet(\Delta)}L_D,
\end{equation}
where the summation on the
right-hand side denotes the subspace spanned by $\{L_D\}_{D\in\facet(\Delta)}$.
Note that if $E\in \Delta$, then $L_E \subset L_D$ for some
$D\in \facet(\Delta)$.  Therefore the right-hand is spanned by 
$L_E$, $E\in \Delta$.

\subsection{The group of invariance of a toric ideal}
\label{subsec:group-of-invariance}

Now we give a definition of the group of invariance of a toric ideal.

Let $S_\cI$ denote the symmetric group on $\cI$, i.e.\ an element
$g\in S_\cI$ is a permutation of the cells of $\cI$. Then
%$S_\coldim$, $\coldim=|\cI|$, acts on the set of cells $\cI$ and then 
$g\in S_\cI$  acts (from the left) on the $|\cI|$-dimensional
rational vector space $\Q^{|\cI|}=\{(y(\Bi))_{\Bi\in \cI}\}$
by the permutation of components: $(gy)(\Bi)=y(g^{-1}(\Bi))$. 
Similarly $g$ acts on the set of indeterminates $\{u_\Bi\}_{\Bi\in \cI}$.
If we regard $g$ as a linear map from $\Q^{|\cI|}$ to itself, then 
it is represented by a permutation matrix.
We denote the permutation matrix also by $g$.
Note that $g$ is orthogonal.
For a given subspace $L\subset \Q^{|\cI|}$, let 
$G_L = \{ g\in S_{|\cI|} \mid g L = L\}$  denote the set-wise stabilizer of $L$.

Let $A$ be a $\rowdim\times \coldim$ rational matrix as in the previous
subsection.  The symmetric group $S_\coldim$ acts on the set of
columns of $A$ and on $\Q^\coldim$.
In \cite{largest} we defined the {\em group of
  invariance} for $A$ as the set-wise stabilizer $G_{\ker A}\subset
S_\coldim$ of $\ker A$. From the viewpoint of toric ideal, the group
of invariance is the set of permutations of the indeterminates, which
leaves the toric ideal invariant. Let $r(A) \subset \Q^\coldim$ denote
the linear space spanned by the rows of $A$.  By Proposition 1 of
\cite{largest}, we have $G_{\ker A}=G_{r(A)}$.  

Our objective is to
understand $G_{\ker A_\Delta}= G_{r(A_\Delta)}$ of a hierarchical
model $\Delta$.

\subsection{A perturbation lemma}

Here we present the following lemma.
% \marginpar{\tiny () or \{\}\\ notation? \\
%  () for a vector\\
%   \{\} for a set}

\begin{lemma} (Perturbation lemma)\ 
Let $n,b$ be positive integers.  There exist $n$ positive integers
$(Y_l)_{l=1}^n$, such that
\begin{equation}
\label{eq:perturbation}
  \{ -b, -b+1, \dots, b-1, b\}^n \ni (c_l)_{l=1}^n
\mapsto \sum_{l=1}^n c_l Y_l
\end{equation}
is injective. Furthermore we can choose  $n$ vectors 
$Y^{(j)}=(Y_l^{(j)})_{l=1}^n$, $j=1,\dots,n$, such that 
(\ref{eq:perturbation}) is injective for each $j$ and
they constitute a basis of the vector space $\Q^n$.
\end{lemma}

\begin{proof}
Let $Y_l^{(j)}=(2b+j)^{l-1}$, $(l,j\in [n])$.  By the uniqueness of the
base $2b+j$ expression of positive integers, the map
$(c_l)_{l=1}^n \mapsto \sum_{l=1}^n c_l Y_l^{(j)}$ is injective.
Furthermore $Y_l^{(j)}$, $j=1,\dots,n$, are linearly independent
in view of the van der Monde determinant.
\end{proof}

In view of the above lemma, we define a generic contingency table belonging
to $r(A_\Delta)$ for 
a given set of cells $\cI$ and a hierarchical model $\Delta$ 
with $\facet(\Delta)=\{D_1, \dots, D_K\}$.

\begin{definition} \label{def:generic}
For $n=\rowdim=\sum_{k=1}^K |\cI_{D_k}|$ and 
$b=\coldim=|\cI|$ choose $(Y_l)_{l=1}^n$  such that 
(\ref{eq:perturbation}) is injective. Decompose  
$(Y_l)_{l=1}^n$ into subvectors of sizes $|\cI_{D_k}|$, $k=1,\dots,K$, as
$
(Y_l)_{l=1}^n = \big((\theta_{D_1}(\Bi_{D_1}))_{\Bi_{D_1}\in \cI_{D_1}},\dots,
(\theta_{D_K}(\Bi_{D_K}))_{\Bi_{D_K}\in \cI_{D_K}}\big)$
and define 
\[
x(\Bi) = \theta_{D_1}(\Bi_{D_1}) + \dots + \theta_{D_K}(\Bi_{D_K}), \quad \Bi\in \cI.
\]
We call this $x$ a generic element of $r(A_\Delta)$.
\end{definition}

Note that an element $g$ of 
the group of invariance $G_{r(A_\Delta)}$ has 
to map a generic element $x$ of $r(A_{\Delta})$ into $r(A_\Delta)$. This fact helps us to determine
$G_{r(A_\Delta)}$.

\section{Group of invariance of hierarchical models}
\label{sec:group-of-invariance-hierarchical}

In this section we first consider a candidate group for the group of invariance
$G_{\ker A_\Delta}$
%=G_{r(A_\Delta)}$ 
and then present our main theorem, which states that the candidate group
is indeed the group of invariance,
provided that the number of levels
$I_j$, $j\in [m]$, are generic.  

% by the perturbation lemma we prove that
% it is the group of invariance, 
%See Theorem \ref{thm:main}.

For $D\subset [m]$
consider a simplicial complex $\Delta^D$, which consists of all
subsets of $D$, and let 
$
\ker A_{\Delta^D} = L_D^{\bot}= \{ y \mid y^+(\Bi_D)=0, \forall \Bi_D \in \cI_D \}
$.  Then 
$
\ker A_\Delta = \bigcap_{D\in \facet(\Delta)} \ker A_{\Delta^D}
$
by  (\ref{eq:kern1}).
% Note that $M_D=\ker A_{\Delta^D}$ is the kernel of the configuration
% $A_{\Delta^D}$ for the hierarchical model $\Delta^D$.
Let $G_D= G_{\ker A_{\Delta^D}}$ denote the group of
invariance for $\Delta^D$.   Then it is easily seen that
%Then we have the following lemma.
% \begin{lemma}
% \label{lem:candidate-group}
% \begin{equation}
%\label{eq:simplex-intersection}
%\bigcap_{D\in \facet(\Delta)} G_D \subset G_{\ker A_\Delta}.
%\end{equation}
%\end{lemma}
%Proof is obvious and omitted. % (or will be added later).
%
\begin{equation}
\label{eq:simplex-intersection}
\bigcap_{D\in \facet(\Delta)} G_D \subset G_{\ker A_\Delta}.
\end{equation}
Therefore we can take $\bigcap_{D\in \facet(\Delta)} G_D$
as a candidate group for the group of invariance $G_{\ker A_\Delta}$.
As we will present an example of sudoku in Section~\ref{sec:wreath}, in general
the inclusion in (\ref{eq:simplex-intersection}) is strict. 
However if the number of levels
$I_j$, $j\in [m]$, are generic, then 
the inclusion in (\ref{eq:simplex-intersection}) is in fact an equality.

Before stating our main theorem, we prove that $G_D=
G_{\ker A_{\Delta^D}}$
is a wreath product of
symmetric groups.
Let $S_{\cI_D}$ denote the symmetric group acting on the
set of $D$-marginal cells and $S_{\cI_{D^C}}$ 
denote the symmetric group acting on the
set of $D^C$-marginal cells, where $D^C$ is the complement of $D$.
Let $(S_{\cI_{D^C}})^{\cI_D}$ denote the set of all functions
from $\cI_D$ to $S_{\cI_{D^C}}$.
Then the wreath product $S_{\cI_{D^C}}\wreath S_{\cI_D}$
is a set $W=S_{\cI_D}\times (S_{\cI_{D^C}})^{\cI_D}$.
The operation of $W$ as a subgroup of $S_{\cI}$ is defined by its action to $\cI$,
where $g=(h,\tilde{h})\in W$
acts on $\Bi\in\cI$ by
$(g\Bi)_D=h\Bi_D$ and $(g\Bi)_{D^C}=\tilde{h}(\Bi_D)\Bi_{D^C}$.
Then we have the following proposition.

\begin{proposition}
\label{prop:1}
The group of invariance $G_D$ for the hierarchical model $\Delta^D$ is given
by the wreath product $S_{\cI_{D^C}} \wreath S_{\cI_D}$.
\end{proposition}

\begin{proof}
For notational simplicity, we prove the proposition for the case
of $m=2$ and $D=\{1\}$ and write $\Bi$ as $(i,j)$.
We denote 
$S_{\cI_D}$ and $S_{\cI_{D^C}}$ by $S_{I_1}$ and $S_{I_2}$, respectively.
The proof for a general case is totally the same 
by the consideration of a ``pseudofactor'' 
%(see Remark \ref{rem:combined} below).
(see Section \ref{sec:hierarchical-facts} for details on pseudofactors).

First we show that $S_{I_2} \wreath
S_{I_1} \subset G_D$.
Let $x\in r(A_{\Delta^D})=L_D$.  Then $x(i,j)=\theta(i)$ for some $\theta$. Let
$g\in S_{I_2} \wreath S_{I_1}$. Then $g(i,j)=(h(i), \tilde h_i(j))$,
where $h\in S_{I_1}$ and $\tilde h_i\in S_{I_2}$ for each $i\in [I_1]$.  Then
\[
(gx)(i,j)=x(g^{-1}(i,j))=\theta(g^{-1}(i,j)_1)=\theta(h^{-1}(i)),
\]
where the subscript ``1'' in $g^{-1}(i,j)_1$ denotes the first component.
Therefore $gx\in L_D$.

We now show the converse $G_D \subset S_{I_2} \wreath
S_{I_1}$.  In order to show this we assume that $x\in L_D$ is generic, i.e.\ 
$\theta(i)$'s are distinct.
Suppose that  $(gx)(i,j)=\theta(g^{-1}(i,j)_1) \in L_D$.  Then
$g^{-1}(i,j)_1$ does not depend on $j$.  Therefore we can write
$g^{-1}(i,j)=(h(i), \bar h(i,j))$. Since $g$ is a bijection, $h$ is a bijection and 
$j\mapsto \bar h(i,j)$  is a bijection for each $i$.
Therefore $g^{-1}\in
S_{I_2} \wreath
S_{I_1}$.
\end{proof}

Now we state the main theorem of this paper.
\begin{theorem}
\label{thm:main}
Consider a hierarchical model $\Delta$. 
Assume that $|\cI_D|$, $D\in \facet(\Delta)$, are
distinct and $I_j > 2$ except for at most one $j\in [m]$. 
Then the group of invariance $G_{\ker A_\Delta}$ is given
by
\begin{equation}
\label{eq:thm1}
G_{\ker A_\Delta} = \bigcap_{D\in \facet(\Delta)} \big(
 S_{\cI_{D^C}} \wreath S_{\cI_D} \big) .
\end{equation}
\end{theorem}

A proof of this theorem is given in Section \ref{sec:proof} after we
establish several important facts on hierarchical models in Section
\ref{sec:hierarchical-facts}.  As seen from the statement of Theorem
\ref{thm:main}, it seems that the case of two-level factors $I_j=2$
needs a special consideration, although the requirements on the levels
in Theorem \ref{thm:main} may be too restrictive.  We discuss these points again
in Section \ref{sec:discussions}.
We will give some examples of Theorem \ref{thm:main}
in Section~\ref{sec:wreath}
after rewriting the right-hand side of (\ref{eq:thm1}).

\section{Some basic facts on hierarchical models}
\label{sec:hierarchical-facts}

%The rest of this section is devoted to the proof of this theorem.

In this section we establish basic facts on hierarchical models.
In particular we are interested  
in the behavior of a hierarchical model, when a maximal simplex is
deleted from $\facet(\Delta)$.
This is because for our proof  of Theorem \ref{thm:main}
we employ the induction on the number $K=|\facet(\Delta)|$ of maximal
interaction effects in $\Delta$.
% Let $\Delta$ be an abstract simplicial
% complex and 

Let $E\subset [m]$.
We first define ``incremental subspaces'' of $L_E$
by 
\begin{equation}
 \label{eq:incremental}
 N_E=L_E\cap\left(\sum_{j\in E}L_{E\setminus\{j\}}\right)^{\bot}
\end{equation}
if $E\neq\emptyset$, and $N_{\emptyset}=L_{\emptyset}$.
Recall that 
$L_E$ is the linear space of tables depending only on the marginal cell $\Bi_E$
and that $\sum_{j\in E}L_{E\setminus\{j\}}$ is
the subspace spanned by $\{L_{E\setminus\{j\}}\}_{j\in E}$.
The following lemma is easily proved and well known
in statistical analysis of variance (ANOVA).

\begin{lemma} \label{lem:increment}
 Let $E$ and $F$ be subsets of $[m]$. Then
 \begin{itemize}
  \item[(1)] $N_E=L_E\cap(\sum_{j\in E}L_{[m]\setminus\{j\}})^{\bot}$.
  \item[(2)] If $E\neq F$, then $N_E\bot N_F$.
  \item[(3)] $L_E=\sum_{F\subset E}N_F$.
  \item[(4)] For any simplicial complex $\Delta$,
	     $r(A_{\Delta})=\sum_{F\in\Delta}N_F$
	     and $\ker A_{\Delta}=\sum_{F\notin\Delta}N_F$.
  \item[(5)] The orthogonal projection $\pi_{N_E}$ onto $N_E$
	     is given by
	     \[
	      (\pi_{N_E}x)(\Bi) = (\pi_{N_E}x)(\Bi_E)
	     = \sum_{F\subset E}
	     \frac{(-1)^{|E\setminus F|}}{|\cI_{F^C}|}
	     x^+(\Bi_F)
	     \]
	     for all $x\in\mathbb{Q}^{\cI}$.
	     Recall that $|E\setminus F|$ is the cardinality of
	     $E\setminus F$.
 \end{itemize}
\end{lemma}

Let $D\in \facet(\Delta)$ be a maximal simplex.  
As in the beginning of Section~\ref{sec:group-of-invariance-hierarchical}
let $\Delta^D$ denote the simplicial complex 
consisting of all subsets of $D$.  Note that $r(A_{\Delta^D})=L_D$.
Now we define $\Delta_{\setminus D}$ by 
``deleting the maximal interaction effects $D$ from $\facet(\Delta)$", i.e.\ 
by
\[
\facet(\Delta_{\setminus D}) = (\facet(\Delta)) \setminus D.
\]
%In the notation of Definition 2.12 of \cite{kozlov},
%$\Delta_{\setminus D} = {\rm dl}_\Delta(\tau)$, where 
%$\tau=D \setminus \cup_{D'\neq D, D'\in \facet(\Delta)} D'$ is  the set of factors
%belonging only to $D$.  
% We have the following lemma, which is a direct
% consequence of ``ANOVA (analysis of variance) decomposition'' 
% of the hierarchical model.
% Using the fact $r(A_{\Delta_{\setminus D}})
% =(\ker A_{\Delta_{\setminus D}})^\perp$ we prove the following lemma.
We have the following proposition.

\begin{proposition} 
\label{prop:induction}
 Let $D\in \facet(\Delta)$.  Then
\begin{align}
\label{eq:induction}
r(A_\Delta) \cap \ker A_{\Delta_{\setminus D}}=
r(A_{\Delta^D})\cap \ker A_{\Delta_{\setminus D}}=
\sum_{E\in\Delta^D\setminus\Delta_{\setminus D}} N_E.
%r(A_{\Delta^D})\cap \ker A_{\Delta^D\cap\Delta_{\setminus D}}.
%
%&=\{ x\in L_D \mid x^+(\Bi_{D \cap D'})=0, 
%\forall \Bi_{D \cap D'} \in \cI_{D\cap D'}, \ 
%\forall D'\in (\facet(\Delta))\setminus D\}.
\end{align}
\end{proposition}

\begin{proof}
By Lemma~\ref{lem:increment},
we have $r(A_{\Delta})=\sum_{E\in\Delta}N_E$
and $\ker A_{\Delta_{\setminus D}}=\sum_{E\notin\Delta_{\setminus D}}N_E$.
Therefore the equalities follow from
the relation $\Delta\setminus\Delta_{\setminus 
 D}=\Delta^D\setminus\Delta_{\setminus D}$.
% $= \Delta^D\setminus(\Delta^D\cap\Delta_{\setminus D})$.
\end{proof}

We next define a partial difference operator.  For $j\in [m]$ and
$x=(x(\Bi))_{\Bi\in \cI}$ 
define
\[
(\partial_j x)(\Bi)=x(\Bi)-x(i_1,\dots, i_{j-1},1,i_{j+1},\dots, i_m), 
\quad \Bi=(i_1,\dots,i_m).
\]
For $E \subset [m]$ define $\partial_E = \prod_{j\in E} \partial_j$.
Note that for two subsets $D,E\subset [m]$, 
$E\not\subset D$, 
we have 
\begin{equation}
\label{eq:partial-difference}
\partial_E \theta_D =0, \qquad \forall \theta_D\in L_D.
\end{equation}
It is obvious that for any $D,E\subset [m]$ and $\theta_D\in L_D$,
we have $\partial_E\theta_D\in L_D$.

Concerning the partial difference operator $\partial_E$
we have the following proposition.

\begin{proposition}
\label{prop:integration}
For all $E\subset [m]$,
$\ker \partial_E =\sum_{F\not\supset E}N_F$.
\end{proposition}

\begin{proof}
We first show that the subspace $\ker\partial_j$
is equal to $L_{[m]\setminus\{j\}}$.
Let $x\in\ker\partial_j$.
Then 
$x(\Bi)=x(i_1,\ldots,i_{j-1},1,i_{j+1},\ldots,i_m)$
and therefore $x\in L_{[m]\setminus\{j\}}$.
Conversely, if $x\in L_{[m]\setminus\{j\}}$,
then $\partial_jx=0$.
Therefore we see that $\ker\partial_j=L_{[m]\setminus\{j\}}$.
Since the operators $\{\partial_j\}_{j\in [m]}$
are mutually commutable projectors (and therefore simultaneously
diagonalizable), 
we have $\ker\partial_E=\sum_{j\in E}\ker\partial_j$.
Therefore, by using Lemma \ref{lem:increment},
\[
 \ker\partial_E
 =\sum_{j\in E}\ker\partial_j
 =\sum_{j\in E}L_{[m]\setminus\{j\}}
 =\sum_{j\in E}\sum_{F\subset [m]\setminus\{j\}}N_F
 =\sum_{F\not\supset E}N_F.
\]
The last equality comes from the fact
that $F\not\supset E$ if and only if
$F\subset [m]\setminus\{j\}$ for some $j\in E$.
\end{proof}

Combining Lemma \ref{lem:increment} and Proposition
\ref{prop:integration}, we have the following proposition.
We will use the proposition with
$\Delta'=\Delta_{\setminus D}$ in the proof of the main theorem.

\begin{proposition}
\label{prop:induction-integration}
Let $\Delta$ and $\Delta'$ be two simplicial complexes
such that $\Delta\supset\Delta'$.
Then $x\in r(A_{\Delta'})$ if and only if $x\in r(A_{\Delta})$
and $\partial_Ex=0$ for all $E\in\Delta\setminus\Delta'$.
%Let $D\in \facet(\Delta)$ and let
% $x=(x(\Bi))_{\Bi\in \cI}
% \in r(A_\Delta)$.  Then
%$x \in r(A_{\Delta_{\setminus D}})$ if and only if
%$\partial_E x=0$ 
%for all
%$E\in \Delta \setminus \Delta_{\setminus D}$.
\end{proposition}

\begin{proof}
The statement is equivalent to
$r(A_{\Delta'})=r(A_{\Delta})\cap(\cap_{E\in\Delta\setminus\Delta'}
 \ker\partial_E)$.
The left-hand side is $\sum_{F\in\Delta'}N_F$.
The right-hand side is
\[
 \left(\sum_{F\in\Delta}N_F\right)\cap
 \left(\bigcap_{E\in\Delta\setminus\Delta'}
 \sum_{F\not\supset E}N_F\right)
 = \sum_{F\in\Delta''}N_F,
\]
where $\Delta''=\{F\in\Delta\mid F\not\supset E,\ \forall 
 E\in\Delta\setminus\Delta'\}$.
It is sufficient to prove that $\Delta'=\Delta''$.
Let $F\in\Delta'$.
Clearly $F\in\Delta$.
Now assume that there exists some $E\in\Delta\setminus\Delta'$
such that $F\supset E$.
Then, since $F\in\Delta'$ and $F\supset E$, we have $E\in\Delta'$.
This contradicts to $E\in\Delta\setminus\Delta'$.
Therefore $F\not\supset E$ for any $E\in\Delta\setminus\Delta'$.
Conversely,
suppose that $F\in\Delta$ and $F\not\supset E$ for any
$E\in\Delta\setminus\Delta'$.
Then, since $F\in\Delta$ and $F\notin\Delta\setminus\Delta'$,
we have $F\in\Delta\setminus(\Delta\setminus\Delta')=\Delta'$.
\end{proof}

In the proof of
Proposition \ref{prop:1}, we treated the combination of factors in $D$
as a single factor and the combination of factors in $D^C$ as another
single factor.  This identification is well known in design of
experiments as a ``pseudofactor'' (e.g.\ \cite{monod-bailey-1992}).
As the last topic of this section we fully discuss the notion of a
pseudofactor and a natural partial order induced on the set of
pseudofactors from the hierarchical model.  The resulting poset plays
an essential role in the next section.

For each  $i\in [m]$ let $\overline{\rm st}_\Delta(i)=\{D\in\Delta\mid \{i\}\cup D\in\Delta\}$
denote the closed star of $\{i\}$, or equivalently the cone over the link of vertex $i$ (see e.g.\ \cite[Definition 2.14]{kozlov}). Let
\begin{equation}
\label{eq:redi}
\fst_\Delta(i)
 = \facet(\overline{\rm st}_\Delta(i))
 = \{ D\in\facet(\Delta)\mid i\in D\}
\end{equation}
denote the facets of $\overline{\rm st}_\Delta(i)$.
If $\fst_\Delta(i)=\fst_\Delta(j)$ we say that $i,j$ belong to the
same pseudofactor and denote this as $i\stackrel{\Delta}{\sim}j$.
%$i \stackrel{\Delta}{\sim}j$ can be alternatively defined as
%\[
%\forall D \in \red\Delta,  \quad i\in D \Leftrightarrow j\in D.
%\]
The relation $\stackrel{\Delta}{\sim}$ is an equivalence relation and
$[m]$ is partitioned into disjoint equivalence classes. We call each equivalence class
a pseudofactor.  In the framework of this paper, we can replace a pseudofactor
by a single factor, although we do not do this in this paper.
Let $\cP$ denote the set of pseudofactors, i.e.
$\cP=[m]/\!\stackrel{\Delta}{\sim}$. For $\rho\in \cP$ let
\[
\fst_\Delta(\rho)= \fst_\Delta(i), \qquad  i\in \rho.
\]
Now we introduce a partial order onto $\cP$ by
\[
\rho  \geq  \rho'  \ \Leftrightarrow \ 
\fst_\Delta(\rho) \supset \fst_\Delta(\rho')
\ \Leftrightarrow\ \overline{\rm st}_\Delta(\rho)\supset\overline{\rm st}_\Delta(\rho').
\]
With this partial order $\cP$ becomes a partially ordered set (poset).
We call this poset the ``pseudofactor poset''
induced by the simplicial complex
$\Delta$.

The pseudofactor poset induced by $\Delta$ is related to the intersection poset.
The intersection poset $\cQ$ of $\facet(\Delta)$
is the set of intersections of $\facet(\Delta)$, that is,
$\cQ=\{\cap_{D\in S}D\mid S\subset\facet(\Delta)\}$.
The order of $\cQ$ is the
reverse inclusion order: $\cap_{D\in S}D\leq \cap_{D\in S'}D$
if $S\subset S'$.
We assume $[m]\in\cQ$
just for convenience.
We show that there is an injective homomorphism from $\cP$ into $\cQ$.
In fact, the following lemma holds.

\begin{lemma} \label{lem:homomorphism}
 Let $V(\rho)=\cup_{\rho'\geq \rho}\rho'$.
 Then 
 $V(\rho)=\cap_{D\in\fst_\Delta(\rho)}D$.
 Furthermore $V$ is an injective homomorphism from $\cP$ into $\cQ$.
\end{lemma}

\begin{proof}
 Let $i\in V(\rho)$.
 Then there exists some $\rho'\geq\rho$ such that $i\in\rho'$.
 This means 
 $\fst_\Delta(i)=\fst_\Delta(\rho')\supset\fst_\Delta(\rho)$.
 Therefore
 $i\in\cap_{D\in\fst_\Delta(\rho)}D$.
 The converse is similarly proved.
 Next we prove that $V$ is homomorphic and injective.
 If $\rho'\geq\rho$ then
 $\fst_\Delta(\rho')\supset\fst_\Delta(\rho)$
 and therefore $V(\rho')=\cap_{D\in\fst_\Delta(\rho')}D\geq
 \cap_{D\in\fst_\Delta(\rho)}D=V(\rho)$ from the definition of the order of $\cQ$.
 If $\rho'\neq \rho$, then 
 $V(\rho')=\cup_{\rho''\geq\rho'}\rho''\neq\cup_{\rho''\geq\rho}\rho''=V(\rho)$.
\end{proof}

We remark that $V$ is not surjective in general.
For example, let $m=3$ and
$\facet(\Delta)=\{\{1,2\},\{2,3\},\{1,3\}\}$.
Then $\cP=\{\{1\},\{2\},\{3\}\}$ with a trivial order
(i.e.\ no two distinct elements are comparable)
and $\cQ=\{\emptyset,\{1\},\{2\},\{3\},\{1,2\},\{2,3\},\{1,3\},\{1,2,3\}\}$.
The homomorphism is $V(\{i\})=\{i\}$ for $i\in\{1,2,3\}$.
Thus $V$ is not surjective.
In other words, the poset $\cQ$ has the same amount of information
as $\facet(\Delta)$ because $\facet(\Delta)=\facet(\cQ\setminus\{[m]\})$,
but the poset $\cP$ loses the information as the example shows.
For description of the group of invariance,
we only need the pseudofactor poset rather than the intersection poset.

\section{A proof of the main theorem}
\label{sec:proof}

%(Rewrite this proof $\rightarrow$ Sei kun)
Now we employ induction on $K=|\facet(\Delta)|$.  
The theorem is
true for $K=1$ by Proposition \ref{prop:1}. Therefore assume that
the theorem holds for $K-1$.
Throughout the proof
we choose
$D\in\facet(\Delta)$ such that $|\cI_D|=\min_{F\in\facet(\Delta)}|\cI_F|$.
We consider deleting $D$ from $\facet(\Delta)$.

Let $x=\sum_{F\in \facet(\Delta)} \theta_F$ be a generic element of
$r(A_\Delta)$ (Definition~\ref{def:generic}).
List the values of $\theta_F$ as
$\alpha_{\Bi_F}=\theta_F(\Bi_F)=\theta_F(\Bi)$,
$\Bi_F\in \cI_F$. Then $x(\Bi)$ can be written as
\begin{equation}
\label{eq:alternative-generic}
x(\Bi)=\sum_{F\in \facet(\Delta)} \sum_{\Bj_F\in \cI_F} \chi_{\Bj_F}(\Bi)
  \alpha_{\Bj_F}, \qquad
\chi_{\Bj_F}(\Bi)=\begin{cases} 1 & {\rm if}\  \Bi_F = \Bj_F, \\
                            0& {\rm otherwise}.
                          \end{cases}
\end{equation}
%In view of Lemma \ref{lem:candidate-group} 
In view of (\ref{eq:simplex-intersection})
it suffices to show
that any $g\in G_{r(A_\Delta)}$ belongs to the right-hand side of 
(\ref{eq:thm1}). 
Fix an arbitrary $g\in G_{\ker A_\Delta}=G_{r(A_{\Delta})}$ and let $y=gx$.
Then $y\in r(A_\Delta)$ and $y$ can
be written as
$y=\sum_{F\in \facet(\Delta)}\eta_F$. Note that at this point we do not
have any relation between $\theta_F$'s and $\eta_F$'s.
Fix an arbitrary $E\in \Delta \setminus \Delta_{\setminus D}$
and take the partial difference with respect to $E$. Then
\begin{equation}
\label{eq:partial}
\partial_E y(\Bi) = \partial_E \eta_{D}(\Bi)
\end{equation}
by (\ref{eq:partial-difference}).  The right-hand side 
$\partial_E \eta_{D}(\Bi)$ depends only  on $\Bi_{D}$. 
The left-hand side $(\partial_E y)(\Bi)$  is a linear combination of $2^{|E|}$ $y(\Bj)$'s 
with the coefficient $1$ for $2^{|E|-1}$ terms and
$-1$ for other $2^{|E|-1}$ terms. Now
$y(\Bj)=x(g^{-1}(\Bj))=(x\circ g^{-1})(\Bj)$.  
We substitute $x(g^{-1}(\Bj))$ by the
right-hand side of (\ref{eq:alternative-generic})  and take the linear
combination.  Then $(\partial_E y)(\Bi)$ is written as
\begin{equation}
\label{eq:partial-E-y}
(\partial_E y)(\Bi)=\sum_{F\in \facet(\Delta)} \sum_{\Bj_F \in \cI_F}
Q_{\Bj_F}(\Bi) \alpha_{\Bj_F},
\end{equation}
where
\[
Q_{\Bj_F}(\Bi)= (\partial_E (\chi_{\Bj_F}\circ g^{-1}))(\Bi)
\in \{-2^{|E|-1}, \dots, 2^{|E|-1} \}.
\]
Since we have taken generic $\alpha_{\Bj_F}$'s, by the perturbation lemma, 
$Q_{\Bj_F}(\Bi)$ is uniquely determined by
$(\partial_E y)(\Bi)$
for each $\Bi$ and for each $F\in \facet(\Delta)$ and $\Bj_F$.
However recall  by
(\ref{eq:partial})
that $(\partial_E y)(\Bi)$ only depends on $\Bi_{D}$.
This implies that 
$Q_{\Bj_F}(\Bi)$ also depends only on $\Bi_{D}$ 
for each $\Bj_F$.  More precisely, if we take the $\Bi_{D}$ marginal of
(\ref{eq:partial-E-y}), then we have
\[
(\partial_E y)^+ (\Bi_{D})= |\cI_{D^C}| (\partial_E y) (\Bi)
= \sum_{F\in \facet(\Delta)} \sum_{\Bj_F \in \cI_F}
Q_{\Bj_F}^+(\Bi_{D}) \alpha_{\Bj_F}.
\]
Therefore by uniqueness  we see that
$Q_{\Bj_F}(\Bi)=Q_{\Bj_F}^+(\Bi_{D})/|\cI_{D^C}|$
depends only on $\Bi_{D}$.

Now we claim that $Q_{\Bj_F}(\Bi)=0$ for all $\Bj_F$, $F\neq D$, and
for all $\Bi\in \cI$.  For readability, we state this as a lemma and
give a proof. Recall that $E\in \Delta \setminus \Delta_{\setminus
  D}$ is arbitrarily fixed and the following lemma holds for any
such $E$.

\begin{lemma}
$Q_{\Bj_F}(\Bi)=0$ for all $\Bj_F\in\cI_F$, $F\in(\facet(\Delta))\setminus\{D\}$, and
for all $\Bi\in \cI$.
\end{lemma}

\begin{proof}
Suppose that there exists some $\Bi^0$ and 
some $\Bj_F$, such that $Q_{\Bj_F}(\Bi^0)\neq 0$. Then, because
$Q_{\Bj_F}(\Bi^0)$ only depends on $\Bi^0_{D}$, 
for this $\Bj_F$ we have
\[
| \{ \Bi \mid Q_{\Bj_F}(\Bi)\neq 0 \}|  \;
\ge \; |\{ \Bi \mid \Bi_{D}=\Bi^0_{D} \}|
=
|\cI_{D^C}| = \frac{|\cI|}{|\cI_{D}|}.
\]
Write 
\[
\cI_{|\Bi^0_{D}}=\{ \Bi \mid \Bi_{D}=\Bi^0_{D} \}
=\{ (\Bi_{D^C}, \Bi_{D}^0)\}_{\Bi_{D^C}\in \cI_{D^C}}.
\]
The $D^C$-component
$\Bi_{D^C}$ of the elements of $\cI_{|\Bi^0_{D}}$ are all distinct.

For $\Bi\in \cI_{|\Bi^0_{D}}$, consider
$Q_{\Bj_F}(\Bi_{D^C}, \Bi_{D}^0)=
(\partial_E (\chi_{\Bj_F}\circ g^{-1}))(\Bi_{D^C}, \Bi_{D}^0)$, 
which is a sum of 
$2^{|E|}$  terms
of the form $\pm \chi_{\Bj_F}(g^{-1}(\Bi'))$.
Since  the operator $\partial_E$ only touches
indices $i_j, j\in E\subset D$, we note that
these terms  $\chi_{\Bj_F}(g^{-1}(\Bi'))$
have the common index $\Bi_{D^C}$, i.e., $\Bi'_{D^C}=\Bi_{D^C}$.
Therefore we can write
\begin{equation}
\label{eq:disjoint}
Q_{\Bj_F}(\Bi_{D^C}, \Bi_{D}^0)= \sum_{\Bj_{D}\in \cI_{D}}
\beta_{\Bj_{D}} \chi_{\Bj_F}(g^{-1}(\Bi_{D^C}, \Bj_{D})),
\end{equation}
where $\beta_{\Bj_{D}}\in \{-1,0,1\}$.
It is important to note
that the sets of cells $\{g^{-1}(\Bi_{D^C},
\Bj_{D})\}_{\Bj_{D}\in \cI_{D}}$ are mutually disjoint for
different values of $\Bi_{D^C}$, because $g^{-1}$ is a bijection on
$\cI$.

Now if $Q_{\Bj_F}(\Bi_{D^C}, \Bi_{D}^0)\neq 0$, 
there exists at least one non-zero term on the right-hand side of
(\ref{eq:disjoint}).  Therefore for each $\Bi_{D^C}$ there exists
$\Bj_{D}$ such that $\chi_{\Bj_F}(g^{-1}(\Bi_{D^C}, \Bj_{D}))=1$.
By the disjointness noted above, it follows that
\[
| \{ \Bi' \mid  \chi_{\Bj_F}(g^{-1}(\Bi'))=1\}| \; \ge \; |\cI_{D^C}|.
\]
On the other hand, by the definition of $\chi_{\Bj_F}$, we have
\[
| \{ \Bi' \mid  \chi_{\Bj_F}(g^{-1}(\Bi'))=1\}|=
| \{ \Bi' \mid  \chi_{\Bj_F}(\Bi')=1\}|=|\cI_{F^C}|.
\]
Combining the above results we have
\[
\frac{|\cI|}{|\cI_F|}=|\cI_{F^C}| \ge |\cI_{D^C}|=
\frac{|\cI|}{|\cI_{D}|} \qquad {\rm or} \qquad
|\cI_{F}| \le  |\cI_{D}|
\]
However we have assumed that $|\cI_{D}|$ is the (unique) minimum among
$|\cI_F|$, $F\in \facet(\Delta)$.  Therefore $F=D$.
\end{proof}

From the above lemma, we have
\begin{equation}
\label{eq:Dk-only}
(\partial_E y)(\Bi)=\sum_{\Bj_{D} \in \cI_{D}}
Q_{\Bj_{D}}(\Bi) \alpha_{\Bj_{D}}.
\end{equation}
We have shown (\ref{eq:Dk-only}) for generic $x$. However, since
(\ref{eq:Dk-only}) is an algebraic relation and all generic tables span 
$r(A_{\Delta})$ by the perturbation lemma,
(\ref{eq:Dk-only}) holds for all $x\in r(A_{\Delta})$.
Now in  (\ref{eq:alternative-generic}) set
$\alpha_{\Bj_{D}}=0$, $\forall \Bj_{D}\in {\cI_{D}}$.
Namely let $x=\sum_{F\in \facet(\Delta), F\neq D} \theta_F$
be any element of $r(A_{\Delta_{\setminus D}})$.
Then  $\partial_E y=\partial_E (x\circ g^{-1})=0$ for
all $E\in \Delta \setminus \Delta_{\setminus
  D}$.
Therefore $y \in r(A_{\Delta_{\setminus D}})$ by Proposition
\ref{prop:induction-integration}. This means that
$g\in G_{r(A_\Delta)}$ has to map every $x\in 
r(A_{\Delta_{\setminus D}})$ into  $r(A_{\Delta_{\setminus D}})$.
In other words, $g\in G_{r(A_{\Delta_{\setminus D}})}$.
By induction assumption we have shown 
\[
g\in \bigcap_{F\in \facet(\Delta), F\neq D} \big(
 S_{\cI_{F^C}} \wreath S_{\cI_F} \big).
\]

Now it remains to show that $g\in S_{\cI_{D^C}} \wreath S_{\cI_{D}}$.
By assumption $g$ maps $r(A_\Delta)$ into itself.  We have shown
that $g$ maps $r(A_{\Delta_{\setminus D}})$ into itself.
Since $g$  is orthogonal as a linear map, it follows that $g$ maps
the subspace $M=r(A_{\Delta})\cap r(A_{\Delta_{\setminus D}})^{\perp}$
into itself.
By Proposition \ref{prop:induction}, we obtain
\[
M=r(A_\Delta)\cap r(A_{\Delta_{\setminus D}})^\perp
=r(A_{\Delta^D})\cap \ker A_{\Delta_{\setminus D}}
=\sum_{E\in\Delta^D\setminus\Delta_{\setminus D}} N_E.
\]
Recall that $N_E$ is the incremental subspace defined by (\ref{eq:incremental}).
Note that $N_D\subset M\subset r(A_{\Delta^D})$.
We claim that
there exists
a table $\phi_D$ in $M$
such that $\phi_D(\Bi_{D})$, $\Bi_{D}\in\cI_{D}$,
are all distinct.
We state this as a lemma and give a proof.

\begin{lemma} \label{lem:distinct}
There exists a table $\phi_D$ in $M$
such that $\phi_D(\Bi_{D})$, $\Bi_{D}\in\cI_{D}$,
are all distinct.
\end{lemma}

\begin{proof}
Consider a generic element $\theta_D$ of $L_D$.
Let $\pi_{N_D}$ denote the orthogonal projection to $N_D$
and put $\phi_D=\pi_{N_D}\theta_D$.
By Lemma~\ref{lem:increment}, the following expression for $\phi_D(\Bi_D)$ holds.
\[
 \phi_D(\Bi_D)
 = \sum_{E\subset D}(-1)^{|D\setminus E|}
 \frac{1}{|\cI_{E^C}|}\theta_{D}^+(\Bi_E).
\]
Recall that $\theta_D^+(\Bi_E)=\sum_{\Bj\in\cI,\Bj_E=\Bi_E}\theta_D(\Bj_D)$.
Multiplying  each side by $|\cI|$, we have
\begin{eqnarray*}
 |\cI| \phi_D(\Bi_D)
 &=& \sum_{E\subset D}(-1)^{|D\setminus E|}
 |\cI_E|\theta_{D}^+(\Bi_E)
 \\
 &=& \sum_{\Bj_D\in\cI_D}C(\Bi_D,\Bj_D)\theta_D(\Bj_D),
\end{eqnarray*}
where
\begin{eqnarray*}
 C(\Bi_D,\Bj_D)
 &=& \sum_{E\subset \mathrm{eq}(\Bi_D,\Bj_D)}(-1)^{|D\setminus E|}|\cI_E||\cI_{D^C}|
 \\
 &=& |\cI_{D^C}|(-1)^{|D\setminus \mathrm{eq}(\Bi_D,\Bj_D)|}
  \prod_{j\in \mathrm{eq}(\Bi_D,\Bj_D)}(I_j-1)
\end{eqnarray*}
and $\mathrm{eq}(\Bi_D,\Bj_D)=\{j\in D\mid i_j=j_j\}$.
Note that $C(\Bi_D,\Bj_D)\in\{-|\cI|,\ldots,|\cI|\}$.
For given $\Bi_D$ and $\Bi'_D$,
if $C(\Bi_D,\Bj_D)\neq C(\Bi'_D,\Bj_D)$ for some $\Bj_D$,
then $\phi_D(\Bi_D)\neq \phi_D(\Bi'_D)$
because $\theta_D$ is generic.
Therefore it is sufficient to prove that
if $\Bi_D\neq\Bi'_D$, then
there exists some $\Bj_D\in\cI_D$ such that 
$C(\Bi_D,\Bj_D)\neq C(\Bi'_D,\Bj_D)$.
Since $I_j$ is greater than 2 except for at most one $j\in [m]$,
we can show that $C(\Bi_D,\Bj_D)=|\cI_{D^C}|\prod_{j\in D}(I_j-1)$
if and only if $\Bi_D=\Bj_D$.
Thus $C(\Bi_D,\Bi_D)\neq C(\Bi'_D,\Bi_D)$ whenever $\Bi_D\neq\Bi'_D$.
This proves the lemma.
\end{proof}

We have proved that
there exists $\phi_D\in N_D\subset M$
such that $\phi_D(\Bi_D)$, $\Bi_D\in\cI_D$,
are all distinct.
Since $g\phi_D\in M\subset r(A_{\Delta^D})$,
the same proof as in Proposition \ref{prop:1} shows that
$g\in S_{\cI_{D^C}} \wreath S_{\cI_{D}}$. 

This completes the proof of Theorem \ref{thm:main}.

\section{The wreath product indexed by the pseudofactor poset} 
\label{sec:wreath}

Although (\ref{eq:thm1}) gives a form of the group of invariance, it is
not yet sufficiently explicit to write down the group of invariance 
for a given hierarchical model.  We can employ the notion of
a wreath product of a partially ordered set of actions to describe 
the group of invariance more explicitly.
The notion of a wreath product of a partially ordered set of actions
has been defined by many authors (\cite{holland-1969},
\cite{wells1976}, \cite{silcock-1977}, \cite{bailey-1983}).
We follow a succinct definition in Section 7 of \cite{wells1976}.
% gives the 
% notion of a wreath product of an partially ordered set of actions.

The poset we use is the pseudofactor poset $(\cP,\leq)$ defined in
Section~\ref{sec:hierarchical-facts}.
Recall that $\cP$ is a partition of $[m]$
and each class $\rho\in\cP$ has $\fst_\Delta(\rho)=\fst_\Delta(i)=\{D\in\facet(\Delta)\mid i\in D\}$,
for $i\in\rho$.
The order relation $\rho\leq\rho'$ on $\cP$
is defined by $\fst_\Delta(\rho)\subset\fst_\Delta(\rho')$.
Recall that $V(\rho)=\cup_{\rho'\geq\rho}\rho'$.
We also define the ancestor set of $\rho$ by
\[
 A(\rho)=\cup_{\rho'>\rho}\rho'=V(\rho)\setminus \rho.
\]
If $A(\rho)=\emptyset$, then we let $\cI_{A(\rho)}$
be a 1-element set, say $\{1\}$.

\begin{definition}[\cite{wells1976}]
The wreath product of the symmetric groups $(S_{\cI_{\rho}})_{\rho\in\cP}$
indexed by the poset $\cP$
is defined by $W=\prod_{\rho\in\cP}(S_{\cI_{\rho}})^{\cI_{A(\rho)}}$,
where $(S_{\cI_{\rho}})^{\cI_{A(\rho)}}$
is the set of all functions from $\cI_{A(\rho)}$ to $S_{\cI_{\rho}}$.
The action of $w=(w_{\rho})_{\rho\in\cP}\in W$ on $\cI$ is defined by
\[
 (w\Bi)_{\rho} = w_{\rho}(\Bi_{A(\rho)})\Bi_{\rho}.
\]
\end{definition}

In the above definition, 
we use the parentheses for evaluating functions (such as $w_{\rho}(\Bi_{A(\rho)})$)
and do not use them for action (such as $w\Bi$).

For example, if $\facet(\Delta)=\{D\}$ and
$\emptyset\subsetneq D\subsetneq [m]$,
then $\cP=\{D,D^C\}$ with the order relation $D>D^C$.
In this case, the wreath product of $(S_{\cI_{\rho}})_{\rho\in\cP}$
is the usual wreath product
$S_{\cI_{D^C}}\wreath S_{\cI_D}$
because $\cI_{A(D)}=\{1\}$ and $\cI_{A(D^C)}=\cI_D$.

The following lemma by \cite{bailey-1983} is useful.

\begin{lemma}[Theorem B of \cite{bailey-1983}] \label{lem:bailey}
The wreath product is characterized as follows.
\[
 \prod_{\rho\in\cP}(S_{\cI_{\rho}})^{\cI_{A(\rho)}}\ =\ \left\{
 g\in S_{\cI}\mid
 (g\Bi)_{V(\rho)}\ \mbox{depends\ only\ on}\ \Bi_{V(\rho)}
 \ \mbox{for\ any}\ \rho\in\cP
 \right\}.
\]
\end{lemma}

The proof of the following lemma is easy and omitted.

\begin{lemma} \label{lem:depends-only}
 Let $A$ and $B$ be two subsets of $[m]$.
 Let $g\in S_{\cI}$.
 Assume that $(g\Bi)_A$ depends only on $\Bi_A$
 and that $(g\Bi)_B$ depends only on $\Bi_B$.
 Then $(g\Bi)_{A\cap B}$ depends only on $\Bi_{A\cap B}$,
 and $(g\Bi)_{A\cup B}$ depends only on $\Bi_{A\cup B}$.
\end{lemma}

Now we establish the following theorem.

\begin{theorem} \label{thm:wreath}
The group of invariance coincides with the wreath product of $(S_{\cI_{\rho}})_{\rho\in\cP}$,
that is,
\begin{equation}
 \bigcap_{D\in\facet(\Delta)}(S_{\cI_{D^C}}\mathop{\rm wr}S_{\cI_D})
 \ =\ \prod_{\rho\in\cP}(S_{\cI_{\rho}})^{\cI_{A(\rho)}}.
 \label{eq:thm-wreath}
\end{equation}
\end{theorem}

\begin{proof}
 By Lemma~\ref{lem:bailey}, the left-hand side in (\ref{eq:thm-wreath}) is equal to
 \[
  \{g\in S_{\cI}\mid
 (g\Bi)_D\ \mbox{depends\ only\ on}\ \Bi_D\ \mbox{for\ any}\ D\in\facet(\Delta)\}.
 \]
 On the other hand, also by Lemma~\ref{lem:bailey}, the right-hand side in (\ref{eq:thm-wreath}) is 
 equal to
 \[
  \{g\in S_{\cI}\mid
 (g\Bi)_{V(\rho)}\ \mbox{depends\ only\ on}\ \Bi_{V(\rho)}
 \ \mbox{for\ any}\ \rho\in\cP\}.
 \]
 Now the equality (\ref{eq:thm-wreath}) is clear
 if one uses Lemma~\ref{lem:depends-only} with two relations
 \[
  D\ =\ \bigcup_{\rho\in\cP, \rho\subset D}V(\rho)\quad\mbox{and}\quad
 V(\rho)\ =\ \bigcap_{D\in\fst_\Delta(\rho)}D.
 \]
 The former one is from the construction of $\cP$.
 The latter one is Lemma~\ref{lem:homomorphism}.
\end{proof}

\begin{corollary} \label{cor:direct-product}
The group of invariance is equal to the direct product
of the symmetric groups $(S_{\cI_{\rho}})_{\rho\in\cP}$
if and only if the poset $\cP$ has the trivial order,
i.e.\ no two distinct elements of $\cP$ are comparable.
\end{corollary}

Let us present some examples.
Here we abbreviate $(S_{\cI_{\rho}})^{\cI_{A(\rho)}}$ to $S_{\rho|A(\rho)}^*$,
and $S_{\cI_{\rho}}$ to $S_{\rho}^*$, respectively.

\begin{example}
 Let $\facet(\Delta)=\{\{1\},\ldots,\{m\}\}$.
 Then $\cP=\{\{1\},\ldots,\{m\}\}$ with the trivial order.
 The wreath product is the direct product
 $W = \prod_{j=1}^m S_{\{j\}}^*$.
\end{example}

\begin{example}
 Let $m=3$ and $\facet(\Delta)=\{\{1\},\{2,3\}\}$.
 In this case, $\{2,3\}$ is a pseudofactor but not a single factor.
 Then the pseudofactor poset is $\{\{1\},\{2,3\}\}$
 with the trivial order.
 The wreath product is $W = S_{\{1\}}^*\times S_{\{2,3\}}^*$.
\end{example}

\begin{example}
 Let $m=3$ and $\facet(\Delta)=\{\{1,2\},\{2,3\}\}$.
 Then the pseudofactor poset is $\{\{1\},\{2\},\{3\}\}$
 with the order relations $\{1\}<\{2\}$ and $\{3\}<\{2\}$
 (no other relations).
 The wreath product is
 $W=S_{\{1\}|\{2\}}^*\times S_{\{2\}}^*\times S_{\{3\}|\{2\}}^*$.
\end{example}

\begin{example}
 Let $m=3$ and $\facet(\Delta)=\{\{1\},\{2\}\}$.
 Note that the factor $\{3\}$ does not appear explicitly.
 Then the pseudofactor poset is $\{\{1\},\{2\},\{3\}\}$
 with the order relations $\{3\}<\{1\}$ and $\{3\}<\{2\}$.
 The wreath product is
 $W=S_{\{1\}}^*\times S_{\{2\}}^*\times S_{\{3\}|\{1,2\}}^*$.
\end{example}

\begin{example}
 Let $m\geq 3$ and $\facet(\Delta)=\{\{1,2\},\{2,3\},\ldots,\{m-1,m\},\{m,1\}\}$.
 Then $\cP=\{\{1\},\ldots,\{m\}\}$ with the trivial order.
 The wreath product is $W=\prod_{j=1}^mS_{\{j\}}^*$.
\end{example}

\begin{example} \label{example:Markov}
 We give an example in that 
 the number of orbits 
 in the minimal Markov basis is smaller 
 under the action of the group of invariance 
 than under the action of the direct-product group.
 Let $m=5$ and $\facet(\Delta)=\{\{1,3\},\{2,4\},\{3,4,5\}\}$.
 This model is decomposable (see \cite{lauritzen1996} for the definition)
 and the minimal vertex separators of the corresponding chordal graph
 are $\{3\}$ and $\{5\}$.
 The pseudofactor poset is
 $\cP=\{\{1\},\{2\},\{3\},\{4\},\{5\}\}$
 with the order relations $\{1\}<\{3\}$, $\{2\}<\{4\}$, $\{5\}<\{3\}$, $\{5\}<\{4\}$
 (no other relations)
 and the wreath product is  given as
 $W=S_{\{1\}|\{3\}}^*\times S_{\{2\}|\{4\}}^*\times S_{\{3\}}^*\times S_{\{4\}}^*\times S_{\{5\}|\{3,4\}}^*$.
 Consider the following two moves
 \begin{eqnarray*}
   M_1 &=& (11111)(12211)-(12111)(11211),
    \\
   M_2 &=& (11112)(12211)-(12112)(11211),
 \end{eqnarray*}
 where the notation follows one in \cite{largest}.
 The moves $M_1$ and $M_2$ are indispensable because they connect
 the following two-element fibers, respectively.
 \begin{eqnarray*}
  \mathcal{F}_1 &=& \{(11111)(12211),(12111)(11211)\},
    \\
  \mathcal{F}_2 &=& \{(11112)(12211),(12112)(11211)\}.
 \end{eqnarray*}
 The direct-product group $\prod_{j=1}^5 S_{I_j}$
 cannot map $M_1$ to $M_2$
 because the number of the distinct levels of the 5th factor in $M_1$
 is different from that in $M_2$.
 On the other hand, the group of invariance
 maps $M_1$ to $M_2$ with a permutation of the 5th factor when 
 the level of the 3rd factor is 1.  
 
 Note that the $S_{\{5\}|\{3,4\}}^*$ involves the set $\{3,4\}$, which is 
 not a maximal clique nor a minimal vertex separator of the chordal graph.
 This shows that the group of invariance can not be described by
 usual notions of decomposition of a chordal graph.

\end{example}

\begin{example}
 Let $m=6$ and $\facet(\Delta)=\{\{1,4,5\},\{2,5,6\},\{3,4,6\}\}$.
 Then the pseudofactor poset is
 $\cP=\{\{1\},\{2\},\{3\},\{4\},\{5\},\{6\}\}$
 with the order relations $\{1\}<\{5\}$, $\{1\}<\{4\}$,
 $\{2\}<\{5\}$, $\{2\}<\{6\}$, $\{3\}<\{4\}$ and $\{3\}<\{6\}$
 (no other relations).
 The wreath product is
 $W=S_{\{1\}|\{4,5\}}^*\times S_{\{2\}|\{5,6\}}^*\times S_{\{3\}|\{4,6\}}^*
 \times S_{\{4\}}^*\times S_{\{5\}}^*\times S_{\{6\}}^*$.
\end{example}

The last example is a counter-example to
the conjecture in the discussion of \cite[Section 5]{largest}.
In our terminology, the conjecture is stated as
``If all pseudofactors are single, i.e.\ $\cP=\{\{1\},\ldots,\{m\}\}$,
and the intersection of $\facet(\Delta)$ is empty,
then the group of invariance is the direct product of the symmetric 
groups on each factor''.
The conjecture is justified if we impose an additional condition that
$\cP$ has the trivial order (see Corollary~\ref{cor:direct-product}).

We show an example in that the inclusion (\ref{eq:simplex-intersection}) is strict.

\begin{example}[Sudoku] \label{ex:sudoku}
 The solution of sudoku is a $9\times 9$ table
 whose each row, each column and each $3\times 3$ block
 contains the digits from 1 to 9 exactly once.
 Following the terminology of \cite{russel-2006},
 we call a ``row'' of $3$ blocks a band
 and a ``column'' of $3$ blocks a stack.
 The solution is considered as a $3\times 3\times 3\times 3\times 9$
 contingency table $x(i,j,k,l,c)$
 where we define $x(i,j,k,l,c)=1$
 if the number $c\in [9]$ is put
 on the $j$-th row of the $i$-th band and the $l$-th column of the $k$-th stack
 and $x(i,j,k,l,c)=0$ otherwise.
 Then the restriction is given by four equations
 \[
  x(i,j,+,+,c)=1,\quad x(+,+,k,l,c)=1,\quad x(i,+,k,+,c)=1,\quad
  x(i,j,k,l,+)=1, 
 \]
 where ``$+$'' denotes taking marginal (sum) over the index.
 The maximal simplices of this model is given by
 \[
  \facet(\Delta)=\{\{1,2,5\},\{3,4,5\},\{1,3,5\},\{1,2,3,4\}\}.
 \]
 The pseudofactor poset is
 $\cP=\{\{1\},\{2\},\{3\},\{4\},\{5\}\}$
 with the order $\{1\}>\{2\}$ and $\{3\}>\{4\}$ (no other relations).
 The wreath product is given by
 \[
  W = S_{\{1\}}^*\times S_{\{2\}|\{1\}}^*\times S_{\{3\}}^*\times 
 S_{\{4\}|\{3\}}^*\times S_{\{5\}}^*,
 \]
 which consists of permutation of bands, permutation of rows in each band,
 permutation of stacks, permutation of columns in each stack
 and permutation of numbers.
 However, the group of invariance $G_{\ker A_{\Delta}}$
 has an additional permutation $f$ defined by
 $f(i,j,k,l,c)=(k,l,i,j,c)$.
 The permutation $f$ does not belong to the 
 wreath product $W$.
 Note that the model does not satisfy the assumption of
 Theorem~\ref{thm:main}
 because $|\cI_{\{1,2,5\}}|=|\cI_{\{3,4,5\}}|=|\cI_{\{1,3,5\}}|
 =|\cI_{\{1,2,3,4\}}|=81$.
 The group generated by $W$ and $f$
 is used to count the number of essentially different solutions
 of sudoku in \cite{russel-2006}.
%\Begin{center}
% \begin{tabular}{|c|c|c||c|c|c||c|c|c|}
%  \hline
%  3& 1& 6& 7& 8& 4& 2& 9& 5\\\hline
%  4& 2& 9& 3& 6& 5& 7& 1& 8\\\hline
%  5& 7& 8& 1& 2& 9& 3& 4& 6\\\hline\hline
%  9& 3& 7& 6& 5& 1& 8& 2& 4\\\hline
%  6& 8& 2& 9& 4& 7& 5& 3& 1\\\hline
%  1& 4& 5& 8& 3& 2& 9& 6& 7\\\hline\hline
%  8& 9& 1& 4& 7& 3& 6& 5& 2\\\hline
%  2& 6& 4& 5& 9& 8& 1& 7& 3\\\hline
%  7& 5& 3& 2& 1& 6& 4& 8& 9\\\hline
% \end{tabular}
%\end{center}
\end{example}

\section{Discussions}
\label{sec:discussions}

We derived an explicit formula of the group of invariance
provided that the number of levels $I_j$, $j\in [m]$, are generic.
In our future work we intend to generalize this result
by weakening the restriction on the number of levels.
We conjecture that under mild regularity conditions the group of invariance is
generated by the wreath product of this paper 
and the permutation of factors with a common number of levels.
However, it seems to be difficult to solve this problem.
For example, as described in Example 3 of \cite{largest},
the group of invariance for
the $2\times 2\times 2$ contingency tables with fixed two-dimensional marginals
is different from the new conjectured candidate group.
In the example, as was pointed out by a referee to
\cite{largest}, the group of invariance is not faithful.
Here an action $G$ to $L$ is called faithful
if the kernel $\{g\in G\mid gx=x,\ \forall x\in L\}$ of the action
consists only of the unit element.
On the other hand, 
we can prove that the group of invariance is faithful
under the assumption of Theorem~\ref{thm:main}.
Indeed, in a similar way to the proof of Lemma~\ref{lem:distinct},
we can show that there exists  a table $\phi\in N_{[m]}\subset \ker A_{\Delta}$
such that $\{\phi(\Bi)\}_{\Bi\in\cI}$ are all distinct.
Therefore if $g\phi=\phi$, then $g$ has to be the identity map.

Random sampling from the group of invariance
is important for performing the Markov Chain Monte Carlo
(MCMC) method on contingency tables.
See \cite{aoki-takemura-2008aism} for details.
In Theorem~\ref{thm:wreath} we rewrote the group of invariance
from an intersection form to a wreath-product form.
The wreath product is useful for random sampling.
Let us briefly describe it.
The wreath product is given by $W=\prod_{\rho\in \cP} (S_{\cI_{\rho}})^{\cI_{A(\rho)}}$.
We show an algorithm to obtain a uniformly random sample $w=(w_{\rho})_{\rho\in\cP}$
from $W$.
Let us number $\cP$ as $\cP=\{\rho_1,\ldots,\rho_l\}$
such that $i<j$ whenever $\rho_i<\rho_j$.
Then, from $i=l$ down to $1$,
we independently generate $w_{\rho_i}(\Bi_{A(\rho_i)})$ from $S_{\cI_{\rho_i}}$
for each $\Bi_{A(\rho_i)}\in\cI_{A(\rho_i)}$.
The resulting element $w=(w_{\rho})_{\rho\in\cP}$ is a uniformly random sample from $W$.
Remark that the intersection form in Theorem~\ref{thm:main}
does not give such a procedure.

\section*{Acknowledgement}

The authors thank the referee for many helpful suggestions,
and in particular, the construction of Example 6.

\bibliographystyle{elsart-num-sort}
\bibliography{invariant}
\end{document}